\documentclass[]{article}

\usepackage{epstopdf}
\usepackage[caption=false]{subfig}
\usepackage[natbibapa,nodoi]{apacite}
\makeatletter
\def\NAT@def@citea{\def\@citea{\NAT@separator}}
\makeatother
\usepackage{amsfonts}
\usepackage{amsthm}
\usepackage{amssymb}
\usepackage{amsmath}
\newtheorem{theorem}{Theorem}

\newtheorem{lemma}[theorem]{Lemma}
\newtheorem{corollary}[theorem]{Corollary}

\newtheorem{remark}[theorem]{Remark}

\usepackage{xparse}
\usepackage[utf8]{inputenc}
\usepackage[english]{babel}
\usepackage{natbib}

\usepackage[nottoc]{tocbibind}
\usepackage[a4paper,bindingoffset=0.2in,%
left=0.2in,right=0.5in,top=1in,bottom=1in,%
footskip=.25in]{geometry}
\usepackage{changepage}
\usepackage[nopar]{lipsum}

\newcommand{\Ebld}{{\mbox{\bf E}}}
\newcommand{\var}{\mbox{\bf Var}}
\usepackage{comment}

\AtBeginDocument{
	\let\oldref\ref
	\renewcommand{\ref}[1]{
		\IfBeginWith{#1}{fig:}%
		{{\color{blue}Figure~\oldref{#1}}}%
		\IfBeginWith{#1}{eqn:}{eq:}%
		{{\color{blue}\oldref{#1}}}%
		{\IfBeginWith{#1}{tab:}{{\color{blue}Table~\oldref{#1}}}{}}}%
}

\usepackage{color}

\usepackage{caption}

\usepackage{mathtools}
\usepackage{xcolor}
\usepackage{comment}
\usepackage{graphicx}
\usepackage{multirow}
\usepackage{adjustbox}
\usepackage[pagewise, displaymath, mathlines]{lineno}
\usepackage{subfig}
\usepackage{authblk}
\providecommand{\keywords}[1]{\textbf{\textit{Keywords:}} #1}

\newtheorem{corollary}[theorem]{Corollary}{\bfseries}{\itshape}
\newtheorem{lemma}[theorem]{Lemma}{\bfseries}{\itshape}
\newtheorem{remark}[theorem]{Remark}{\bfseries}{\itshape}

\title{Estimating Linear Mixed Effects Models with Truncated Normally Distributed Random Effects}
\author[1]{Hao Chen\thanks{Corresponding Author: hao.x.chen@nielseniq.com}}
\author[1]{Lanshan Han\thanks{lanshan.han@nielseniq.com}}
\author[1,2]{Alvin Lim\thanks{aclim2@emory.edu}}
\affil[1]{Retail Product Research \& Development, NielsenIQ, Chicago, IL 60606}
\affil[2]{Emory University Goizueta Business School, Atlanta, GA 30322}
\date{}  

\begin{document}
\maketitle				
\begin{abstract} \label{sec:abstract}
It is proved that the sum of $ n $ independent but non-identically distributed doubly truncated Normal distributions converges in distribution to a Normal distribution. It is also shown how the result can be applied in estimating a constrained mixed effects model.
\end{abstract}
		
\noindent \keywords{Truncated Normal Distribution; Lindeberg-Feller Theorem; Lindeberg Condition; Constrained Mixed Effects Model}		
		
\section{Motivation}
It is our observation that modern statistical models are heavily dependent on Normal distributions. For example, consider a simple linear regression model:
\begin{equation}
y_i = \beta_0 + \beta_1 x_i + \epsilon_i,
\end{equation}
where the error term $ \epsilon_i $ is assumed to follow $ N(0, \sigma^2) $, so it follows that $ y_i \sim N(\beta_0 + \beta_1 x_i, \sigma^2) $. Hence, maximum likelihood estimation (MLE) can be applied to estimate the unknown parameters $ \beta_0, \beta_1, \sigma^2 $. Another example is the linear mixed effects (LME) model (\cite{mcculloch2014generalized}), say we have a mixed intercept and slope model given below
\begin{equation}
y_{i, j} = \left(\beta_0 + \beta_{0, j} \right)  + (\beta_1 + \beta_{1, j}) x_{i, j} + \epsilon_{i, j},
\end{equation}
where $ i $ and $ j $ indicate the row number and group number, respectively. It is assumed that
\begin{equation}
\beta_{0, j} \sim N(0, \eta_0^2), \beta_{1, j} \sim N(0, \eta_1^2) , \epsilon_{i, j} \sim N(0, \sigma^2)
\end{equation}
and $  \beta_{0, j}, \beta_{1, j}, \epsilon_{i, j} $ are independent to each other. In other words, only with the Normality assumptions of both the error term and the random effects can we have the result that $ y_{i, j} $ is also Normally distributed, upon which nearly all of the modern statistical inference methods are built. The underlying reason is simple: suppose $ x_i  \sim N(\mu_i, \sigma_i^2)$, it is straightforward to show the weighted sum of independent $  \alpha_i x_i  $ is still Normally distributed with $ \alpha_i $ as known constant. Mathematically,
\begin{equation}
\sum_{i=1}^{n} \alpha_i x_i \sim N\left(\sum_{i=1}^{n} \alpha_i \mu_i, \sum_{i=1}^{n} \alpha_i^2 \sigma_i^2  \right).
\end{equation}
The above property is rather neat and elegant. Without such a well-behaved property, the analytical expression of the exact distribution of $ y_{i, j} $ will be unavailable under the LME model. In addition, many other continuous distributions do not have such a nice behavior. That being said, a Normal distribution is unconstrained. Technically speaking, under the LME model with Normality assumptions, there is no control on the sign of the overall slope $ \beta_1 + \beta_{1, j} $ even if researchers have some prior knowledge about its range. \\

We propose to use the truncated Normal distribution (\cite{horrace2005some}) on the random effects so that the overall slope $ \beta_1 + \beta_{1, j} $ will be bounded. However, significant difficulty has been observed by switching from the Normality assumption to the truncated Normality assumption: the sum of truncated Normal distribution is analytically intractable, and there is a lack of study on its large sample property. Hence, in this note we attempt to show that the sum of $ n $ independent but non-identically distributed doubly truncated Normal (DTN) distributions converges in distribution to a Normal distribution. Therefore, inference based on the Normality of $ y_{i, j}$ can still be applied when $ n $ is sufficiently large. \\

\cite{horrace2005some} studied one-sided truncated Normal (TN) distribution, and the authors presented some analytical results about it. \cite{robert1995simulation} talked about how to simulate truncated Normal variables. More recently, \cite{cha2015re} discussed more properties about TN in his PhD thesis. The rest of this article is organized as follows. We will give more analytical results about TN and DTN in Section \oldref{sec:pre}. The main results are presented in Section \oldref{sec:main}, and how the results can be applied to estimate an constrained LME model is discussed in Section \oldref{sec:app}.

\section{Preliminaries} \label{sec:pre}
A truncated Normal (TN) distribution is parameterized by $ 4 $ parameters: location, $ \mu $; scale, $ \eta $; lower bound $ a $; upper bound $ b $. The Normal distribution is a special case of it when $ a = -\infty  $ and $  b = \infty $. The probability density function (PDF) of a $\mathcal {TN}(\mu, \eta^2, [a,b])$, with $\eta>0$, is given by
$$f_{\mathcal {TN}}(x; \mu,\eta^2, a, b) \, = \, \left\{
\begin{array}{ll}
\displaystyle{ \frac{1}{\eta} \frac{\phi\left(\xi\right)}{\Phi\left(b^{\prime}\right)-\Phi\left(a^{\prime}\right)} }, & x\in[a,b]\\[5pt]
0, & \mbox{otherwise}
\end{array}\right.,$$
where $\phi(\cdot)$ and $\Phi(\cdot)$ are PDF and cumulative distribution function (CDF) of the standard Normal distribution, i.e.,
$$\phi(\xi) \, = \, \frac{1}{\sqrt{2\pi}} \exp\left(-\frac{1}{2}\xi^2\right), \,\, \mbox{and} \,\, \Phi(\xi) \, = \, \frac{1}{2} \left[1+\mbox{erf}\left(\frac{\xi}{\sqrt 2}\right)\right],$$
respectively, and
$$\xi \, \triangleq \, \frac{x-\mu}{\eta}, \, \, a^{\prime} \, \triangleq \, \frac{a-\mu}{\eta},\, \, \mbox{ and } \,\,b^{\prime} \, \triangleq \, \frac{b-\mu}{\eta}.$$
The mean and variance of $x \sim \mathcal{TN}(\mu, \eta^2, [a,b])$ are known and given by \cite{olive2008applied}:
$$\Ebld[x] \, = \, \mu + \frac{\phi(a^{\prime}) - \phi(b^{\prime})}{\Phi(b^{\prime}) - \Phi(a^{\prime})} \eta, $$
$$\var[x] \, = \, \eta^2 \left[ 1+ \frac{a^{\prime} \phi(a^{\prime}) - b^{\prime} \phi(b^{\prime})}{\Phi(b^{\prime}) - \Phi(a^{\prime})} -\left(\frac{\phi(a^{\prime})-\phi(b^{\prime})}{\Phi(b^{\prime}) - \Phi(a^{\prime})}\right)^2\right]. $$
In this paper, we are particularly interested in a special case, namely the symmetric doubly truncated normal (DTN) distributions $\mathcal {TN} (\mu, \eta^2, [\mu - \rho \eta, \mu+\rho \eta] )$ with $\rho>0$, denoted by $\mathcal {DTN}(\mu, \eta^2, \rho)$. It is a special case of a TN with $ a = \mu - \rho \eta, b = \mu + \rho \eta $, i.e., the lower bound and upper bound are symmetric around mean $ \mu $. The properties of a DTN distrbution is given by Lemma \oldref{lm:DTN_properties1}.

\begin{lemma}\label{lm:DTN_properties1}
	Suppose $x \sim \mathcal{DTN}(\mu, \eta^2, \rho)$, the following results hold
	\begin{itemize}
		\item[(i).] The density function is $$f_{\mathcal{DTN}}(x;\mu, \eta^2, \rho) \,= \, \left\{ \begin{array}{ll}\displaystyle{\frac{1}{\eta} \frac{\phi(\xi)}{2\Phi(\rho)-1} }, & x\in[\mu-\rho \eta, \mu+\rho \eta]\\
		0, & \mbox{otherwise} \end{array}\right\}$$
		\item[(ii).] The expectation is  $$\Ebld[x] \, = \, \mu, $$
		\item[(iii).] The variance is $$ \var[x] \, = \, \eta^2 \left[1 - \frac{2\rho \phi(\rho)}{2\Phi(\rho)-1}\right],   $$
	\end{itemize}
\end{lemma}

The proof is omitted as it is straightforward to verify the above results. Note that we define DTN distributions with $\rho>0$. In fact, when $\rho=0$, it becomes a deterministic value, and hence the variance is $0$. This is indeed consistent with the fact that
$$ \lim_{\rho\rightarrow 0} \left[ 1 - \frac{2\rho \phi(\rho)}{2\Phi(\rho)-1} \right] \, = \, 1 - \lim_{\rho \rightarrow 0 }\frac{2\rho \phi(\rho)}{2\Phi(\rho)-1} \, = \, 1- \lim_{\rho \rightarrow 0} \frac{2\phi(\rho) + 2\rho \phi'(\rho)}{2\phi(\rho)} \, = \, 0$$
where the second equal sign is due to L'H$\hat{o}$pital's rule.
We report some properties regarding the DTN distribution in Lemma \oldref{lm:DTN_properties2}.

\begin{lemma}\label{lm:DTN_properties2}
	Let $x \sim \mathcal{DTN}(\mu, \eta^2, \rho)$ with $\rho>0, \eta > 0$, the following results hold.
	\begin{itemize}
		\item[(i).] $x-\mu \sim \mathcal{DTN}(0, \eta^2, \rho)$.
		\item[(ii).] For any $x, y \in[\mu-\rho\eta, \mu+\rho \eta]$, if $x+y = 2 \mu$ then $f_{\mathcal{DTN}}(x;\mu, \eta^2, \rho) = f_{\mathcal{DTN}}(y;\mu, \eta^2, \rho)$.
		\item[(iii).] $\var[x] \leq \eta^2$.
		\item[(iv).] Suppose $x' \sim \mathcal{DTN}(\mu, \eta^2, \rho')$, then $\var[x] \leq \var[x']$ if $\rho \leq \rho'$.
			\item [(v).] If $x \sim \mathcal{DTN}(0, \eta^2, \rho)$. Define $ x^{\prime} = k_0 + k_1 x $ with $ k_1 \ne 0$, where $ k_0, k_1 $ are finite real numbers, then $  x^{\prime} \sim \mathcal{DTN}(k_0, k_1^2\eta^2, \rho) $.
	\end{itemize}
\end{lemma}
\proof
(i) and (ii) are obvious from properties of standard Normal distribution PDF $\phi(\xi)$. To prove (iii), we notice that $2\Phi(\rho) > 1$ for all $\rho>0$. Since $\phi(\rho)>0$, we have
$$\frac{2\rho \phi(\rho)}{2\Phi(\rho)-1} > 0, \,\, \forall \rho>0. $$ Therefore,
$$ \var[x] \, = \, \eta^2 \left[1 - \frac{2\rho \phi(\rho)}{2\Phi(\rho)-1}\right] \, \leq \, \eta^2.$$
To prove (iv), we examine the following function
$$ g(\rho) \, \triangleq \, \frac{2\rho \phi(\rho)}{2\Phi(\rho)-1}.$$
It is clear this function is differentiable on $(0,+\infty)$. Noticing that the derivate of $\phi(\rho)$, $\phi'(\rho) = - \rho \phi(\rho)$, the derivative of $g(\rho)$ is written as
\begin{eqnarray}
g'(\rho) & = & \frac{2\phi(\rho) + 2\rho \phi'(\rho)}{2\Phi(\rho)-1} - \frac{4\rho [\phi(\rho)]^2}{[2\Phi(\rho)-1]^2} \nonumber\\
& = & \frac{[2\phi(\rho) + 2\rho \phi'(\rho)][2\Phi(\rho)-1]-4\rho [\phi(\rho)]^2}{[2\Phi(\rho)-1]^2} \nonumber\\
& = & \frac{2\phi(\rho)\left[ (1-\rho^2)(2 \Phi(\rho)-1) - 2\rho \phi(\rho)\right]}{[2\Phi(\rho)-1]^2} \nonumber
\end{eqnarray}
We further let
$$t(\rho) \, \triangleq \, (1-\rho^2)(2 \Phi(\rho)-1) - 2\rho \phi(\rho).$$
It is clear that $t(\rho)$ is continuous in $\rho$ and $t(0)=0$. We also have
$$t'(\rho) \, = \, -2\rho(2\Phi(\rho)-1) < 0,$$
for all $\rho>0$. Therefore, $t(\rho)<0$ for all $\rho>0$. It implies that $g'(\rho)<0$ for all $\rho>0$. Therefore, $g(\rho)$ is a monotonically decreasing function on $(0,\infty)$. Hence (iv) holds readily. For (v), it is straightforward to verify that the PDF of $x'$ is given by
\begin{eqnarray}
f(x^{\prime}) &=&
\left\{
\begin{array}{ll}
\frac{1}{k_1} f\left(\frac{x - k_0}{k_1}\right) & \mbox{if } k_1 > 0 \\[5pt]
-\frac{1}{k_1} f\left(\frac{x - k_0}{k_1}\right) & \mbox{if } k_1 < 0 \\[5pt]
\end{array}
\right. \nonumber\\
&=& \frac{1}{|k_1|} \frac{1}{\eta} \frac{\phi(  \frac{x-k_0 - 0}{k_1\eta})}{2\Phi(\rho) - 1} \nonumber \\
&=& \frac{1}{|k_1| \eta} \frac{\phi(\frac{x-k_0}{k_1\eta})}{2\Phi(\rho) - 1} \nonumber
\end{eqnarray}
The last equation is the PDF of $ \mathcal{DTN}(k_0, k_1^2\eta^2, \rho)  $. \qedhere

\section{Main Results} \label{sec:main}
It is worth pointing out that, while the sum of independent non-identically distributed Normal random variables is Normally distributed, it is not the case for DTNs. The exact distribution of the sum of independent non-identically DTNs is analytically intractable. However, the following Normality results hold.

\begin{theorem} \label{th:clt}
	For every $ x_i \sim \mathcal{DTN}(\mu_i, \eta_i^2, \rho_i)$, the random variables making up the collection $  \mathbf{X}_n = \{ x_i:  1 \le i \le n \} $ are independent with the following conditions.
	\begin{itemize}
		\item $ \mu_i $ are finite, i.e., $ \max_{1 \le i \le n} \mu_i < +\infty$
		\item $ \rho_i $ is bounded from below by $ \underline{\rho} > 0 $
		\item $ \eta_i $ is bounded from below and above by $ \underline{\eta} > 0 $ and $ \bar{\eta} < + \infty $, respectively.
 	\end{itemize}
	Then
	$$ \frac{1}{t_n} \sum_{i=1}^{n} \left(x_i - \mu_i\right) \, \xrightarrow{d} \, \mathcal N(0, 1), $$
	as $n\rightarrow \infty$,
	where $$t_n^2 = \sum_{i=1}^n \var[x_i].$$
\end{theorem}

\proof
For the proof, we will use the well known Lindeberg-Feller theorem \citep{zolotarev1967generalization}:
\textit{
	Suppose that $x_1,x_2,\cdots$ are independent random variables such that $\Ebld[x_i] =\mu_i$ and $\var[x_i] = \sigma_i^2 < +\infty$ for all $i=1,2,\cdots$. Define:
	\begin{eqnarray}
	y_i & = & x_i - \mu_i, \nonumber \\
	s_n^2 & = & \sum_{i=1}^n \var\left[y_i\right] \, = \, \sum_{i=1}^n \sigma_i^2. \nonumber
	\end{eqnarray}
	If the Lindeberg condition
	\begin{equation}\label{eq:linderberg}
	\mbox{for every }\epsilon>0, \frac{1}{s_n^2} \sum_{i=1}^n \Ebld\left[y_i^2 \cdot \mathbf 1_{|y_i|\geq \epsilon s_n}\right] \, \rightarrow \, 0 \,\,\mbox{as } n\rightarrow \infty
	\end{equation}
	is satisfied, then $$\frac{\sum_{i=1}^{n}\left(x_i - \mu_i\right)}{s_n} \, \xrightarrow{d} \, \mathcal N(0,1).$$ 
}
For Theorem~\oldref{th:clt} to hold, it suffices to verify the Lindeberg condition (\oldref{eq:linderberg}). First, by Lemma~\oldref{lm:DTN_properties2}, item (iii), we have
$$\var[y_i]\leq \eta_i^2 < \infty.$$
Next, since $\rho_i \geq \underline{\rho}$ for each $i=1,2,\cdots$, by Lemma~\oldref{lm:DTN_properties2} item (iv), we have
$$\var[y_i] \, \geq \, \eta_i^2 \left[1 - \frac{2\underline{\rho}\phi(\underline{\rho})}{2\Phi(\underline{\rho})-1}\right] \, \geq \, \underline\eta^2 \left[1 - \frac{2\underline{\rho} \phi(\underline{\rho})}{2\Phi(\underline{\rho})-1}\right] \, \triangleq \, v^2,$$
where $  v = \underline{\eta} \sqrt{1 - \frac{2\underline{\rho} \phi(\underline{\rho})}{2\Phi(\underline{\rho})-1}} > 0 $. It follows that $s_n^2 \geq n v^2$ for all $n=1,2,\cdots$. By Lemma \oldref{lm:DTN_properties2},
$y_i \, \sim \, \mathcal{DTN}(0,\eta_i,\rho_i)$. Therefore, with $ u_i \triangleq \frac{y_i}{\eta_i}$, for any given $\epsilon>0$ and for each $i=1,2,\cdots$, we have
\begin{eqnarray}
\Ebld\left[ y_i^2 \cdot \mathbf 1_{ |y_i|>\epsilon s_n }\right] & = & \int_{-\rho_i \eta_i}^{\rho_i \eta_i} y_i^2 f_{\mathcal{DTN}}\left(y_i;0, \eta_i, \rho_i\right) \cdot \mathbf 1_{ |y_i|>\epsilon s_n } d y_i \nonumber \\
& = & \left\{
\begin{array}{ll}
0 & \mbox{if } \epsilon s_n \geq \rho_i\eta_i \\[5pt]
2 \displaystyle {\int_{\epsilon s_n}^{\rho_i\eta_i} y_i^2 f_{\mathcal{DTN}}\left(y_i;0, \eta_i, \rho_i\right) d y_i }& \mbox{if } \epsilon s_n < \rho_i\eta_i
\end{array}
\right. \nonumber\\
& \leq & 2 \int_{\epsilon s_n}^{\infty}  y_i^2 \frac{1}{\eta_i (2 \Phi(\rho_i)-1)} \phi\left(\frac{y_i}{\eta_i}\right) d{y_i}\nonumber\\
&=& \frac{2}{\eta_i} \frac{1}{2\Phi(\rho_i) - 1} \frac{1}{\sqrt{2\pi}}
\int_{\epsilon s_n}^{\infty} {y_i}^2 \exp{\left(-\frac{y_i^2}{2\eta_i^2}\right)} d y_i \nonumber \\
&=& \frac{2}{\eta_i} \frac{1}{2\Phi(\rho_i) - 1} \frac{1}{\sqrt{2\pi}}
\left(\eta_i^3
\int_{\frac{\epsilon s_n}{\eta_i}}^{\infty} u_i^2 \exp\left(-\frac{u_i^2}{2}\right) du_i\right) \nonumber \\
&=& \eta_i^2 \frac{1}{2\Phi(\rho_i) - 1} \sqrt{\frac{2}{\pi}}
\int_{\frac{\epsilon s_n}{\eta_i}}^{\infty} u_i^2 \exp\left(-\frac{u_i^2}{2}\right) du_i \nonumber \\
&\leq& \overline \eta^2 \frac{1}{2\Phi(\underline \rho) - 1} \sqrt{\frac{2}{\pi}}
\int_{\frac{\epsilon \sqrt n v}{\overline \eta}}^{\infty} u_i^2 \exp{\left(-\frac{u_i^2}{2}\right)} du_i. \nonumber
\end{eqnarray}
Notice that
$$  \int u_i^2 \exp{\left(-\frac{u_i^2}{2}\right)} du_i = \sqrt{\frac{\pi}{2}} \text{erf}\left(\frac{u_i}{\sqrt{2}}\right) - u_i\exp{\left(-\frac{u_i^2}{2}\right)}.     $$
Hence, we have
\begin{eqnarray}
\lim_{n\rightarrow \infty}\frac{1}{s_n^2} \sum_{i=1}^n \Ebld[y_i^2 \cdot \mathbf 1_{|y_i|\geq \epsilon s_n}] & \leq & \lim_{n\rightarrow \infty}\frac{1}{s_n^2} \sum_{i=1}^n \left[\overline \eta^2 \frac{1}{2\Phi(\underline \rho) - 1} \sqrt{\frac{2}{\pi}}
\int_{\frac{\epsilon \sqrt n v}{\overline \eta}}^{\infty} u_i^2 \exp\left(-\frac{u_i^2}{2}\right) d u_i \right] \nonumber\\
& \leq & \lim_{n\rightarrow \infty} \frac{\overline \eta^2 }{n v^2} \frac{1}{2\Phi(\underline \rho) - 1} \sqrt{\frac{2}{\pi}} \sum_{i=1}^n \left[ \sqrt{\frac{\pi}{2}} \left(1-\mbox{erf}\left(\frac{\epsilon \sqrt n v}{\overline \eta}\right) \right) + \frac{\epsilon \sqrt n v}{\overline \eta} \exp\left(-\frac{\epsilon^2 n v^2}{2 \overline \eta^2}\right) \right]\nonumber \\
& = & \lim_{n\rightarrow \infty} \frac{\overline \eta^2 }{v^2} \frac{1}{2\Phi(\underline \rho) - 1} \sqrt{\frac{2}{\pi}} \left[ \sqrt{\frac{\pi}{2}} \left(1-\mbox{erf}\left(\frac{\epsilon \sqrt n v}{\overline \eta}\right) \right) + \frac{\epsilon \sqrt n v}{\overline \eta} \exp\left(-\frac{\epsilon^2 n v^2}{2 \overline \eta^2}\right) \right] \nonumber \\
& = & 0, \nonumber
\end{eqnarray}
where the last equality is due to the fact that since $ \epsilon > 0 $ is finite, and $ v, \overline{\eta} $ are fixed constants,
$$  \lim\limits_{n \rightarrow \infty} \frac{\epsilon \sqrt{n} v}{\overline{\eta}} \rightarrow \infty $$
$$\lim_{n\rightarrow \infty} \mbox{erf}\left(\frac{\epsilon \sqrt n v}{\overline \eta}\right) \, = \, 1 \, \, \mbox{and} \,\, \lim_{z\rightarrow \infty } z \exp\left(-\frac{z^2}{2}\right) \, = \, 0. $$
\qed \\

Moreover, it is also straightforward to verify the following corollary to Theorem~\oldref{th:clt}.
\begin{corollary}\label{co:weighted_CLT}
	Let $x_i \sim \mathcal{DTN}(\mu_i, \eta_i^2, \rho_i)$, $i=1,2,\cdots$ be independent with $\mu_i$'s, $\eta_i$'s, and $\rho_i$'s satisfying conditions in Theorem \oldref{th:clt}. Let $\beta_i, i=1,2,\cdots,$ be real numbers and the absolute values are bounded from below and above, i.e., there exist $\overline \beta$ and $\underline \beta$ satisfying $ 0 < \underline{\beta} \leq |\beta_i| \leq \overline{\beta} < +\infty$ for all $i=1,2,\cdots$. Then,
	$$ \frac{1}{t_n} \sum_{i=1}^{n} \beta_i (x_i - \mu_i) \, \xrightarrow{d} \, \mathcal N(0, 1), $$
	as $n\rightarrow \infty$,  where $$ t_n^2 = \sum_{i=1}^n \beta_i^2 \var[x_i].$$ \qed
\end{corollary}
Corollary~\oldref{co:weighted_CLT} indicates that the (weighted) sum of finitely many independent but non-identically distributed DTNs converges in distribution to a Normal distribution.


\section{Application to Constrained Mixed Effects Model} \label{sec:app}
Suppose there are $ g $ groups, indexed by $\ell = 1, \ldots, g$, the mixed effects model (\cite{mcculloch2014generalized}) is given by
\begin{equation}\label{eq:mixed_reformulation_matrix_proposed}
\boldsymbol y^\ell \, =\, X^\ell \boldsymbol \beta + Z^\ell \boldsymbol \gamma^\ell + \boldsymbol \varepsilon^\ell,
\end{equation}
where
\begin{equation}
\boldsymbol \varepsilon^\ell \sim \mathcal{N}( \mathbf{0}_{m_\ell}, \sigma^2 \mathbf{I}_{m_\ell})
\end{equation}
and $ m_\ell $ is the sample size for group $ \ell $, the total size is $ m = \sum_{\ell=1}^{g} m_{\ell} $, $ \mathbf{0}_{m_\ell} $ is a size $ m_\ell $ column vector with $ 0 $ as all of its elements. $ \mathbf{I}_{m_\ell} $ is a identity matrix with size $ m_\ell $. For the random effects $ \gamma_{\ell, i} $, we assume they are independent and follow the distribution
\begin{equation}
\gamma_{\ell, i} \sim \mathcal{DTN}(0, \varsigma_i^2, [-\beta_i, \beta_i ]), i = 1, \ldots, p,
\end{equation}
where $ \beta_i > 0 $, for each $i=1,\cdots,p$. $ p $ is the number of columns for which the random effects are considered. Each $ \gamma_{\ell, i} $ is mathematically constrained within its corresponding $[-\beta_i, \beta_i ]$. Hence, the overall coefficient of group $ \ell $ and column $ i $ calculated as $ \beta_i + \gamma_{\ell, i} \ge 0  $. This way, we can guarantee that the overall coefficient will be non-negative. One can follow a similar procedure if a non-positive sign is needed. Following the results in Section \oldref{sec:main}, we have
\begin{eqnarray}
\Ebld[y^{\ell}|X^{\ell}, \boldsymbol{\beta}] & = & X^{\ell} \boldsymbol \beta, \nonumber \\
\var[y^{\ell}|X^{\ell}, Z^{\ell}, \boldsymbol{\beta}] & = & Z^{\ell} \Lambda (Z^{\ell})^T + \sigma^2 \mathbf I_{m_{\ell}},\nonumber
\end{eqnarray}
where $$ \Lambda \, = \, \mbox{diag}\left[\left(\varsigma_i^2  \left[1-\frac{2\rho_i \phi(\rho_i)}{2\Phi(\rho_i)-1}\right]\right)_{i=1}^p\right].$$
Therefore, we have
$$  y^{\ell}  \xrightarrow{d} \mathcal N(X^{\ell}\boldsymbol \beta, Z^{\ell} \Lambda (Z^{\ell})^T + \sigma^2 \mathbf I_{m_{\ell}}), $$
and MLE can be used for parameter estimation.

\bibliographystyle{plainnat}
\bibliography{sample}

\end{document}